\documentclass[a4paper,10pt]{amsart}
\usepackage{bbm,amsmath,amsfonts,amssymb,mathrsfs}
\usepackage{yfonts}
\usepackage{dsfont}
\usepackage{version}
\usepackage{cite}

\newcommand{\RingS}[1]{\Complexes(s)[\mathfrak{S}_{#1}]}
\newcommand{\RingG}[1]{k[\mathfrak{G}_{#1}]}
\newcommand{\Ring}{{\mathcal{R}}}
\newcommand{\domain}{\mathbb{X}}

\newcommand{\field}{\ensuremath{k}}
\newcommand{\Naturals}{\ensuremath{\mathbb{N}}}
\newcommand{\Integers}{{\mathbb{Z}}}
\newcommand{\Rationals}{\ensuremath{\mathbb{Q}}}
\newcommand{\Reals}{\ensuremath{\mathbb{R}}}
\newcommand{\Complexes}{\ensuremath{\mathbb{C}}}

\excludeversion{ExtendedVersion}
\newcommand{\Mikusinski}{Mikusi\'nski}
\newcommand{\Bezout}{B\'ezout}

\newcommand{\ie}{\emph{i.e.}}

\newcommand{\seeeg}{\emph{see, e.g.,} }

\newcommand{\eye}{I}
\newcommand{\MikReg}{\mathcal{M}_{\mathrm{R}}}
\newcommand{\MikComp}{\mathcal{M}_0}
\newcommand{\Mik}{\mathcal{M}}
\newcommand{\LaplaceTrans}{\mathscr{L}}
\newcommand{\entire}{\mathcal{O}}
\newcommand{\Boem}{\mathcal{B}}
\newcommand{\DistComp}{\Cinfty'}
\newcommand{\Cinfty}{\mathcal{E}}
\newcommand{\CinftyComp}{\mathcal{D}}
\newcommand{\Dist}{\CinftyComp'}
\newcommand{\Gevrey}[1]{\Cinfty_{#1}}
\newcommand{\DistGevrey}[1]{\Dist_{#1}}

\DeclareMathOperator{\rank}{tr}

\DeclareMathOperator{\Hom}{Hom}
\DeclareMathOperator{\supp}{supp}

\newcommand{\Vect}[1]{ \boldsymbol{#1} }

\newtheorem{definition}{Definition}[section]
\newtheorem{example}{Example}[section]
\newtheorem{notation}{Notation}[section]

\newtheorem{remark}{Remark}[section]
\newtheorem{lemma}{Lemma}[section]
\newtheorem{theorem}{Theorem}[section]
\newtheorem{proposition}{Proposition}[section]
\newtheorem{corollary}{Corollary}[section]
\newtheorem{theoremanddefinition}{Theorem and Definition}[section]
\newtheorem{definitionandproposition}{Definition and Proposition}[section]
\newcommand{\bw}{\Vect{w}}
\newcommand{\bvarK}{ {\ensuremath{\Vect{c}}}} 
\newcommand{\bvarU}{ {\ensuremath{\Vect{u}}}} 
 
\newcommand{\varK}{ {\ensuremath{c}}} 

\begin{document}
\title[Controllability of networks of second order p.d.e.]{Controllability of networks of one-dimensional second order p.d.e. -- an algebraic approach}
\author[F. Woittennek]{Frank Woittennek}
\address{Laboratoire d'Informatique\\
        \'Ecole Polytechnique\\
            91128 Palaiseau Cedex\\ France}
\email{frank@woittennek.de}
\author[H. Mounier]{Hugues Mounier}
\address{D\'epartement AXIS\\
             Institut d'\'Electronique Fondamentale \\
            B\^at. 220, Universit\'e Paris-Sud\\
             91405 Orsay\\
            France}
\email{hugues.mounier@ief.u-psud.fr}

\begin{abstract}
  We discuss controllability of systems that are initially given by boundary
  coupled p.d.e.\ of second order. Those systems may be described by modules
  over a certain subring $\mathcal{R}$ of the ring $\mathcal{M}_0$ of
  \Mikusinski\ operators with compact support.  We show that the ring $\mathcal{R}$ is a
  B\'ezout domain. This property  is utilized in order to derive algebraic 
  and trajectory related controllability results.
\end{abstract}
\maketitle


%
%
%
\section{Introduction}
The solution of control design problems is, in general, preceded by a
controllability analysis of the system under consideration. While for
linear finite dimensional systems, both algebraic and analytic
controllability notions are used in parallel, the analysis of
infinite dimensional systems is dominated by (functional) analytic
methods \cite{CurZwa95}. The latter approach has proven to be
useful, in particular for the analysis of state space controllability,
\ie, the possibility of steering the system under consideration from a
given initial state to a desired final state.  For example,
controllability of the same class of systems as considered in the
present contribution has been analyzed this way in
\cite{LagLeuSchmi94,DagZua2006}.
However, when focusing on trajectory tracking problems, the behavioural
controllability notion due to Willems \cite{Wil91} is an interesting
alternative to classical state space controllability. The connections
between this approach and the algebraic system properties have been
pointed out in the past for several classes of distributed parameter
systems, as for so called multidimensional systems
\cite{PomQua99,Zerz00book}, for delay systems
\cite{Gluesing97siam,Mounier1998fm}, and more general convolutional
systems \cite{VetZam00siam,VetZam02laa}. In addition, the algebraic
viewpoint is closely related to parametrization of trajectories of the
system under consideration. This constructive nature makes the
approach very attractive for applications, in particular for open loop
control design. Finally, paying attention to particular structural
properties of the models under consideration may result in a deeper
understanding of the respective class of systems.

From the algebraic (module theoretic) viewpoint, as used within this
contribution, a linear system is a finitely presented module.  This
notion was first introduced by Fliess  for linear
finite dimensional systems in \cite{Fli90lin}.
For this class of systems, the freeness
of the  module corresponds to the flatness of the system under
consideration in the sense of the theory of nonlinear finite
dimensional systems while its basis corresponds to a flat output
\cite{FLMR95ijc}. Moreover, torsion freeness, \ie,
the absence of autonomous subsystems, is equivalent to freeness.  The
module theoretic approach is applicable to systems with distributed
parameters and lumped controls as well: The convolutional equations
associated with a given boundary value problem serve as defining
relations for the system module, the latter defined over a certain ring
of generalized functions.  Such coefficient rings are generally not
principal ideal domains.  For this reason, the two basic
controllability related properties, torsion freeness and freeness, are
not necessarily equivalent.  An approach to circumvent the problems
caused by this ``lack of structure'' is the concept of $\pi$-freeness,
which relies on localization and was at first developed for linear
delay systems \cite{FliMou98cocv}. This way a basis can be introduced
at least within an appropriate extension of the module under
consideration. The approach has been proven to be very useful for both
trajectory planning and open loop control design
\cite{Woi07diss,RudWoi08ijc,Rud03habil,RWW03md,PetRou01siam,%
  LMR00ijrnc,WoiRud03cocv}. Nevertheless it seems to be difficult to
compare such purely algebraic controllability notions to the
behavioural ones. For this reason, within the present contribution we
do not use localization.

This paper addresses the development of an algebraic approach to the
controllability of networks of spatially one-dimensional parabolic and
hyperbolic constant coefficient p.d.e.\ of second order.  Here, by a
network we understand a system consisting of several branches each of
which is governed by a system of p.d.e.\ and which are coupled via the
boundary conditions.  We investigate, through algebraic properties of
the coefficient rings obtained for considered class of p.d.e., the
related controllabilities of the associated system module and
establish some controllability results including module theoretic and
behavioral ones. In accordance with
\cite{Woi07diss,RudWoi08ijc,RudWoi2007at}, we use the general solution
of the Cauchy problem with respect to space in order to rewrite the
given model as a linear system of convolutional equations. The latter
are regarded as the defining relations of a finitely presented
module. It turns out that the coefficient ring of this module, a
subring of the ring of \Mikusinski\ operators with compact support
introduced in \cite{Boehme1973ams}, is a \Bezout\ domain, \ie, every
finitely generated ideal is principal. An algorithm enabling us to
calculate the generator of a given finitely generated ideal is
presented within this paper.  This latter result is strongly inspired
by those derived in \cite{BreLoi96msca,Gluesing97siam} for particular
rings of distributed delay operators which in our setting may arise
from the wave equation.  The derived properties of the coefficient
ring allow us to decompose the system module into a free module and a
torsion module.  Finally, from these algebraic results, we deduce the
trajectorian controllability of the free submodule in the sense of
\cite{FliMou02ima} and its behavioural controllability in the sense of
\cite{Wil91}.

The paper is organized as follows. In Section \ref{bvp_to_sys} we
introduce the class of models considered as systems of p.d.e.\ which
are coupled via their boundary conditions. We show how to pass from
this model to a system of convolution equations giving rise to our
module theoretic setting.  Section \ref{sec:ring} is devoted to the
study of the coefficient ring of this module.  In Section
\ref{sec:controllability} we obtain several controllability results
for the systems under consideration. Finally, in Section \ref{sec:example},
we apply the method to a system example of two boundary coupled p.d.e.

\section{Boundary value problems as convolutional systems}\label{bvp_to_sys}
\subsection{Models considered}\label{sec:model}
We assume, that the model equations for the distributed variables in
$\Vect{w}_{1},\dots,\Vect{w}_{l}$ and the lumped variables in
$\Vect{u}=(u_1,\dots,u_m)$ are given by
\begin{subequations}\label{eq:bvp}
  \begin{equation}\label{eq:pde}
    \begin{split}
      \partial_x\Vect{w}_i&=A_i\Vect{w}_i+B_i \bvarU,\quad
      \Vect{w}_{i}:\Omega_i\rightarrow \Boem^2,\quad \Vect{u}\in\Boem^m\\
      &\qquad\qquad A_i\in(\Reals[s])^{2\times2},\quad
      B_i\in(\Reals[s])^{2\times m}
    \end{split}
  \end{equation}
  where $\Boem$ denotes an appropriate space of
  Boehmians\footnote{Alternatively, one could use other spaces of
    generalized functions given in the inclusion chain
    $\Gevrey{\rho}\subset\Cinfty\subset\Dist\subset\DistGevrey{\rho}\subset\Boem.$
    Here, as usual, $\Cinfty$ denotes the space of infinitely
    differentiable functions, $\Dist$ are the Schwartz
    distributions. Moreover, $\Gevrey{\rho}$ and $\DistGevrey{\rho}$
    are spaces of Gevrey-functions and Gevrey ultra-distributions of
    order $\rho$, respectively.  However, in order to avoid
    distinctions of cases, we shall use the space of Boehmians.}
  (\seeeg\ \cite{Mikusinski1988amh,BurMikNem2005,Mik00} and
  App.~\ref{sec:miku}) and $s$ is the differentiation operator with
  respect to time.  The assumptions which are crucial for the
  applicability of our approach are twofold.  First, we assume that
  all the matrices $A_1,\dots,A_l$ give raise to the same
  characteristic polynomial, namely,
  \begin{equation}\label{eq:chareq}
    \det(\lambda\eye- A_i)=\lambda^2-\sigma,\quad \sigma=a s^2+ bs+c\ne0,\quad a,b,c\in\Reals,\quad a\geqslant0.
  \end{equation}
  Additionally, we require the intervals $\Omega_1,\dots,\Omega_l$ of
  definition of the above differential equations to be rationally
  dependent. More precisely, we assume the $\Omega_i$ ($i=1,\dots,l$)
  to be given by an open neighbourhood of
  \begin{equation}
    \label{eq:length}
    \tilde\Omega_i=[x_{i,0},x_{i,1}],\quad \ell_i=x_{i,1}-x_{i,0}=q_i\ell, \quad q_i\in\Rationals,\quad \ell\in\Reals.
  \end{equation}
In the following, and without
further loss of generality, we assume $x_{i,0}=0$.
  The model is completed by the boundary conditions
  \begin{equation}
    \sum_{i=0}^l L_{i}\Vect{w}_i(0)+R_{i}\Vect{w}_i(\ell_i)+D \Vect{u}=0\label{eq:bc}
  \end{equation}
\end{subequations}
where $D\in (\Reals[s])^{q\times m} $ and $L_{i},
R_{i}\in(\Reals[s])^{q\times 2}$.  
\begin{remark}
  In a more general setting, instead of the boundary conditions \eqref{eq:bc}, one
  could consider auxiliary conditions of the form
$$
\sum_{i=0}^l Q_{i}(\Vect{w}_i)+D \Vect{u}=0.
$$
Here, 
$$
Q_i(\Vect{w}_i)=\sum_{j=0}^\nu L_{i,j}\Vect{w}_i(\alpha_{i,j}\ell)+
\sum_{j=1}^\mu\int_{\Omega_{i,j}} Q_{i,j}^{\star}(x) \bw_i(x)
dx
$$
with $L_{i,j}\in(\Reals[s])^{q\times 2}$,
$Q_{i,j}^\star\in(\Reals[s,x])^{q\times 2}$,
$\Omega_{i}\supset\Omega_{i,j}=[\beta_{i,j,1}\ell,\beta_{i,j,2}\ell]$,
$\alpha_{i,j},\beta_{i,j,k}\in\Rationals\cap \Omega_i$,
$\mu,\nu\in\Naturals$.
\end{remark}
\subsection{Solution of the inital value problem}\label{sec:inival}
This section deals with the solution of a single initial value problem
of the form (\ref{eq:pde}) with initial conditions given at $x=\xi$,
\ie,
\begin{equation}\label{eq:ivp:general}
 \partial_x\bw=A\bw+B\bvarU, \quad \bw(\xi)=\bw_\xi
\end{equation}
with $A$, $B$ having the same properties as $A_i$, $B_i$
($i=1,\dots,l$) introduced within the previous section. To this end, we
start with the initial value problem
$$
  (\partial_x^2  -\sigma) S(x)=0,\quad  S(x)=0,\quad(\partial_x S)(x)=1,
$$
associated with the characteristic equation (\ref{eq:chareq}).
It is well known that this equation has a unique operational solution
as long as the principal part $\partial_x^2-as^2$ of the differential
operator $\partial_x^2-\sigma$ is hyperbolic w.r.t.\ the parallels of
$x=0$. This was implicitly required above by assuming $a\geqslant
0$ in (\ref{eq:chareq}). Moreover, under these assumptions the operator $S$ as well as its
derivative $C=\partial_xS$ correspond to infinitely differentiable
functions mapping $\Omega$ to the ring $\Mik_0$ of \Mikusinski\
operators with compact support\footnote{More specifically, instead of
  stating that $S(x),C(x)\in\Mik_0$, one could distinguish the cases
  $a>0$ and $a=0$. In the first case both operators, $S(x)$ and
  $C(x)$, correspond to distributions of order zero with compact
  support, while in the latter case they correspond to
  ultra-distributions of Gevrey order 2 and support in $0$. Both of
  these spaces may be embedded in $\Mik_0$.} (see App.~\ref{sec:miku}
and \cite{Boehme1973ams} for results related to the support of
\Mikusinski\ operators and App.~\ref{sec:app:representation} for
explicit expressions for $C(x)$ and $S(x)$).

Using the above defined operational functions one easily verifies that
the (unique) solution $x\mapsto\Phi(x,\xi)$ of the initial value
problem
$$
 \partial_x\Phi(x,\xi)=A \Phi(x,\xi),\quad \Phi(\xi)=\eye,
$$
with $\eye$ denoting the identity,
is given by
\begin{equation}
\label{eq:Phi}
 \Phi(x,\xi)=A S(x-\xi)+\eye C(x-\xi).
\end{equation}
From the uniqueness of the solution one deduces the addition formula
$$\Phi(x,\xi)\Phi(\xi,\zeta)=\Phi(x,\zeta).$$ 
For $A$ the companion matrix of the
characteristic polynomial, \ie,
\begin{equation}\label{eq:companion}
  A=
\begin{pmatrix}
  0&1\\
  \sigma&0
\end{pmatrix},\quad \Phi(x,\xi)=\begin{pmatrix}
  C(x)& S(x)\\
  \sigma S(x)& C(x)
\end{pmatrix},
\end{equation}
this yields in particular
\begin{ExtendedVersion}
  \begin{equation*}
    \begin{pmatrix}
      C(x)& S(x)\\
      \sigma S(x)& C(x)
    \end{pmatrix}
    \begin{pmatrix}
      C(y)& S(y)\\
      \sigma S(y)& C(y)
    \end{pmatrix}
    =\begin{pmatrix}
      C(x+y)& S(x+y)\\
      \sigma S(x+y)& C(x+y).
    \end{pmatrix}
  \end{equation*}
  and thus
\end{ExtendedVersion}
\begin{equation}\label{eq:addition_theorems}
  C(x+y)=C(x)C(y)+\sigma S(x)S(y),\quad
  S(x+y)=C(x)S(y)+ S(x)C(y).
\end{equation}
The solution of the initial value problem associated with the inhomogeneous
equation
\begin{equation}\label{eq:sol:diff}
  \partial_x\Psi(x,\xi)=A\Psi(x,\xi)+B
\end{equation}
with homogeneous initial conditions, prescribed at $x=\xi$, is obtained using
the well known variation of constants method. This yields
\begin{equation}
  \label{eq:Psi}\Psi(x,\xi)=\int_\xi^x \Phi(x,\zeta) d\zeta B.
\end{equation}
Thus, the general solution of the initial value  problem (\ref{eq:ivp:general})
reads
$$
\bw(x)=\Phi(x,\xi)\bw_\xi+\Psi(x,\xi)\bvarU.
$$
\begin{ExtendedVersion}
or equivalentely 
$$
  \bw(x)=W(x,\xi)\bvarK,\quad W(x,\xi)=
  \begin{pmatrix}
    \Phi(x,\xi)&\Psi(x,\psi)
  \end{pmatrix},\quad \bvarK_\xi=
  \begin{pmatrix}
    \bw_\xi\\\bvarU
  \end{pmatrix}.
$$
\end{ExtendedVersion}
The entries of the matrix $\Phi$ belong clearly to
$\Complexes[s,C,S]$.  Contrary, according to \eqref{eq:Psi}, the
entries of $\Psi$ may contain also the integrals of $S$ and $C$.
However, if $\sigma\ne 0$ then $A$ is invertible over $\Complexes(s)$.
Thus, using the fact that $\partial_z \Phi(x,z)=-\Phi(x,z)A$ those
integrals can be expressed as 
$$\int_{\xi}^x \Phi(x,\zeta)d\zeta=-\int_{\xi}^x \partial_\zeta\Phi(x,\zeta)A^{-1}d\zeta=(\Phi(x,\xi)-\eye)A^{-1}.$$
Choosing $A$ as in \eqref{eq:companion} one obtains in particular
\begin{ExtendedVersion}
  \begin{equation*}
    \int_0^x
    \begin{pmatrix}
      C(x-\zeta)& S(x-\zeta)\\
      \sigma S(x-\zeta)& C(x-\zeta)
    \end{pmatrix}d\zeta
=    \int_{0}^x
    \begin{pmatrix}
      C(\zeta)& S(\zeta)\\
      \sigma S(\zeta)& C(\zeta)
    \end{pmatrix}d\zeta
    =\begin{pmatrix}
      S(x)&(C(x)-1)/\sigma\\
      C(x)-1&S(x)
    \end{pmatrix}
  \end{equation*}
  and thus
\end{ExtendedVersion}
\begin{equation*}
 \int_0^x C(\zeta)dx = S(x),\quad \int_0^x S(\zeta)dx = (C(x)-1)/\sigma.
\end{equation*}
Later on, the latter equations will essentially ease our controllability analysis.

\subsection{A module presented by a system of convolution equations}\label{sec:module}
In the previous section we have discussed the solutions of the initial value
problems associated with the equations (\ref{eq:pde}). In the sequel these
results are used in order to define an algebraic structure representing 
the model under consideration.
 
Substituting the general solutions of the initial value problems into the
boundary conditions, one obtains the following linear system of equations:
\begin{equation}
  \label{eq:presentation}
   \bw(x)=W_{\Vect{\xi}}(x)\bvarK_{\Vect{\xi}},\quad P_{\Vect{\xi}}\bvarK_{\Vect{\xi}}=0.
\end{equation}
Here, $\Vect{\xi}=(\xi_1,\dots,\xi_n)$, $\Vect{\varK}_{\Vect{\xi}}^T=(\Vect{w}_{1}^T(\xi_1)\cdots\Vect{w}_{l}^T(\xi_l),\Vect{u}^T)$, 
\begin{align*}
   W_{\Vect{\xi}}&=
 \begin{pmatrix}
   \Phi_1(x,\xi_1)&0&0&\Psi_1(x,\xi_1)\\
    0   &\ddots&0&\vdots\\
   0&\cdots&\Phi_l(x,\xi_l)&\Psi_l(x,\xi_l)
 \end{pmatrix}, \quad P_{\Vect{\xi}}=\big(P_{\Vect{\xi},1}\cdots P_{\Vect{\xi},l+1}\big)
\end{align*}
with 
\begin{align*}
  P_{\Vect{\xi},i}&=L_{i}\Phi_i(0,\xi_i)+R_i\Phi_i(\ell_{i},\xi_i),  \quad i=1,\dots, l\\ 
  P_{\Vect{\xi},l+1}&=D+\sum_{i=1}^lL_i\Psi_i(0,\xi_{i})+R_i\Psi_i(\ell_{i},\xi_{i}).
\end{align*}

A suitable algebraic object for the representation of the model under
consideration should contain all the system variables, \ie, the
distributed variables in $\bw$, their values as well as their (spatial
and temporal) derivatives. Moreover, it should reflect not only the
structure imposed by the original boundary value problem
(\ref{eq:pde}) but also that imposed by the solution of the according
initial value problems, \ie, by equation (\ref{eq:presentation}).  In
order to analyze the model we will, therefore, use a module generated
by the variables collected in $\bvarK_{\xi}$ with the presentation
given in (\ref{eq:presentation})
\cite{FliessMounier04esta,FliMou99,FliMou98cocv,Mou95these}. The
choice of the coefficient ring, which has to contain at least the
entries of $W_{\Vect{\xi}}$ and $P_{\Vect{\xi}}$, is discussed below.

According to the previous section, the entries of $W_{\Vect{\xi}}$ and those of
$P_{\Vect{\xi}}$ are composed of the functions $C$ and $S$ mapping $\Reals$ to
$\MikComp$ and all the values of these functions.  Moreover, the matrices may
involve  values of the spatial integrals of $C$ and $S$, too. Thus a possible choice
for the coefficient ring is the ring $\Ring^I_\Reals[s,\mathfrak{S},\mathfrak{S}^I]$. 
Here, for any $\mathbb{X}\subseteq\Reals$,  $\Ring^I_\mathbb{X}=[\mathfrak{S}_\mathbb{X},\mathfrak{S}^I_\mathbb{X}]$, with
\begin{equation*}
  \begin{alignedat}{2}
     \mathfrak{S}&=\{C,S\},&\mathfrak{S}_{\mathbb{X}}&=\{C(z\ell),S(z\ell)|z \in \mathbb{X}\}, \\
     \mathfrak{S}^I&=\{C^I,S^I\},& \quad\mathfrak{S}^I_{\mathbb{X}}&=\{C^I(z\ell),S^I(z\ell)|z\in \mathbb{X}\},
  \end{alignedat}
\end{equation*}
$\ell$ defined as in \eqref{eq:length}, and
$$S^I(x)=\int_0^x S(\zeta)d\zeta,\quad  C^I(x)=\int_0^x C(\zeta)d\zeta.$$
Inspired by the results given in
\cite{Mounier1998fm,BreLoi96msca,Gluesing97siam}, and in view of the
simplification of the analysis of the module properties, instead of the ring
$\Ring_\mathbb{R}^I$, we will use a slightly larger ring, given by
$\Ring_\Reals=\Complexes(s)[\mathfrak{S}_\Reals]\cap\MikComp$.
\begin{definition}\label{def:Sigma}
  The \emph{convolutional system} $\Sigma$ associated with the boundary value
  problem \eqref{eq:bvp} is the module generated by the elements of
  $\bvarK_{\Vect{\xi}}$ over $\Ring=\Ring_{\Reals}[\mathfrak{S},\mathfrak{S}^I]$
  with presentation matrix $P_{\Vect{\xi}}$.  By $\Sigma_\Reals$
  (resp.\ $\Sigma_\Rationals$) we denote the same system but viewed as a module over
  $\Ring_\Reals$ (resp.\ $\Ring_\Rationals$).
\end{definition}
One easily verifies that $\Sigma$ does not depend on the choice of $\Vect{\xi}$
(cf.\ \cite{Woi07diss,RudWoi08ijc}). In view of the assumed mutual rational
dependence of the lengths $\ell_{1},\dots,\ell_{l}$, for the analysis of the
system properties, we will start with the system $\Sigma_\Rationals$, \ie, a
system containing the values of the distributed variables at rational multiples
of $\ell$ only.  However, having analyzed the properties of $\Sigma_\Rationals$,
we may pass to $\Sigma_\Reals$ (resp.\ $\Sigma$)  by an extension of scalars, \ie,
$\Sigma_\Reals=\Ring_\Reals\otimes_{\Ring_\Rationals}\Sigma_\Rationals$ (resp.\ $\Sigma=\Ring\otimes_{\Ring_\Rationals}\Sigma_\Rationals$).

\section{The operator ring  $\Ring_{\Rationals}$ is a \Bezout\ domain}\label{sec:ring}
In this section we study the structures of the ideals within the rings
$\Ring_\Rationals$. To this end, we first establish some results on the
ideals in $\RingS{\Rationals}$ and $\RingS{\Naturals}$.  

\subsection{Ideals  in  $\RingS{\Rationals}$ and $\RingS{\Naturals}$}\label{sec:ring:field}
This section is devoted to the analysis of the ideals in
$\RingS{\Rationals}$ and $\RingS{\Naturals}$.    The
obtained results will be the key for the analysis of the properties of
$\Ring_{\Rationals}$ and $\Ring_{\Naturals}$ introduced in Section
\ref{sec:module}. 

In the following, we will replace $\Complexes(s)$ by any field
$\field$. Moreover, the generating set $\mathfrak{S}_{\mathbb{X}}$ may
be replaced by any set $\mathfrak{G}_{\domain}$ the elements of which
satisfy the addition formulas derived in Section \ref{sec:inival}.
\begin{definition}
  Let $\field$ be a field and $\domain$ an additive subgroup of $\Reals$.  
  By $\mathfrak{G}_{\domain}$ we denote the set $\{S_a,C_a|a\in \domain \}$
  the elements of which are subject to the following relations ($\sigma\in\field$):
  \begin{subequations}\label{eq:hypfun:propA}
    \begin{equation}
      C_aC_b\pm \sigma S_aS_b=C_{a\pm b},\quad S_aC_b\pm C_aS_b=S_{a \pm b}\label{eq:hypfun:prop:1}
     \end{equation}
     and
    \begin{equation}
      C_{0}=1,\quad S_{0}=0.\label{eq:hypfun:prop:2}
     \end{equation}
\end{subequations}
  \end{definition}

From the above definition one easily deduces 
\begin{subequations}\label{eq:hypfun:propB}
    \begin{equation}\label{eq:hypfun:prop:4}
      C_{a}=C_{-a},\quad S_{a}=-S_{-a}
     \end{equation}
as well as 
     \begin{equation}\label{eq:hypfun:prop:7}
       \begin{gathered}
         2C_aC_b=C_{a+b}+C_{a-b},\: 2\sigma S_aS_b=C_{a+b}-C_{a-b},\:
         2C_aS_b=S_{a+b}-S_{a-b}.
       \end{gathered}
    \end{equation}
\end{subequations}
   Moreover, any element  $r\in\RingG{\domain}$ can be written in the form
   \begin{equation}\label{eq:poly}
     r=\sum_{i=0}^n a_{\alpha_i} C_{\alpha_i} +     b_{\alpha_i} S_{\alpha_i},\quad n\in\Naturals, \quad a_{\alpha_i},b_{\alpha_i}\in \field, \quad \alpha_i\in \domain^+
   \end{equation}
  where $\domain^+=\{|\alpha|:\alpha\in\domain\}$. Finally the units in $\RingG{\domain}$ belong to $\field$.

In the following, it is necessary to distinguish the cases where the
equation $\lambda^2-\sigma=0$ has a solution over $\field$ or not. For
our application this is clearly equivalent to the question whether the
roots of the characteristic equation (\ref{eq:chareq}) belong to
$\Reals[s]$.  The necessity to distinguish these cases is explained by
the following simple example which, in addition, shows that the cases
$\mathbb{X}=\Naturals$ and $\mathbb{X}=\Rationals$ need to be analyzed
separately.
\begin{example}\label{ex:ideals}
  Consider the ideal  $\mathfrak{I}=(a,b)$, $a=S_1$, $b=C_1+1$.  Over $\field[\mathfrak{S}_{\Rationals}]$ we have
  $$a=S_1=2C_{1/2}S_{1/2},\quad b=C_1+1=2C_{1/2}^2.$$  Thus, both generators belong to $(C_{1/2})$ which, conversely,
  belongs to $\mathfrak{I}$ since $2C_{1/2}=-\sigma S_{1/2}
  a+C_{1/2}b$.  The ideal $\mathfrak{I}$ is, therefore, generated by $C_{1/2}$ which
  does not belong to $\field[\mathfrak{S}_{\Naturals}]$ if
  $\lambda^2-\sigma$ is irreducible over $k$.  However, the situation
  is different if $\sqrt{\sigma}$ belongs to $\field$. From the
  relations given in (\ref{eq:hypfun:propA}) and
  (\ref{eq:hypfun:propB}), it follows
  immediately that
  $$
  (C_{1/2}+\sqrt{\sigma}\,S_{1/2})(C_{1/2}- \sqrt{\sigma}\,S_{1/2})=1.
  $$
  Over $\field[\mathfrak{S}_{\Rationals}]$, $C_{1/2}$ can be  factorized as 
  $$
   C_{1/2}=(C_{1/2}+\sqrt{\sigma}\,S_{1/2})(1+C_1-\sqrt{\sigma} S_1)/2.
  $$
  The element $C_{1/2}$ is, thus, 
  associated with $1+C_1-\sqrt{\sigma} S_1$ which indeed belongs to
  $\field[\mathfrak{S}_{\Naturals}]$.
\end{example}

\subsubsection{The polynomial $\lambda^2-\sigma$ is reducible over $k$}\label{sec:reducible}
  \begin{proposition}
    The ring $\RingG{\Naturals}$ is a PID.
  \end{proposition}
  \begin{proof}
    From the addition formulas given in (\ref{eq:hypfun:propA}) and (\ref{eq:hypfun:propB}), 
    it follows  that $\RingG{\Naturals}$ is isomorphic
    to $\field[S_1,C_1]$ which, in turn, is  isomorphic to 
    $\field[z^{-1},z]$ by
    $$
      S_1\mapsto \frac{z^{-1}-z^{1}}{\lambda},\quad
      C_1\mapsto {z^{-1}+z^{1}}.
     $$
     The latter ring is Euclidean  with the norm function given by the difference of the
     degrees of the monomials of maximal and minimal degree w.r.t.\ $z$.
  \end{proof}
 \begin{corollary}\label{corollary:bezout1}
   The ring $\RingG{\Rationals}$ is a \Bezout\ domain.
  \end{corollary}
  \begin{remark}
    Note that, for the same reason as given in Remark~\ref{rem:nopid},
    $\RingG{\Rationals}$ is not a PID.
  \end{remark}
  \begin{remark}
    Having in mind that in our application $\sigma$ is given according
    to equation (\ref{eq:chareq}) and $\field=\Complexes(s)$, for
    $\sqrt{\sigma}\in \field$, say $\sqrt{\sigma}=\lambda$, the
    operators $C(x)$ and $S(x)$ introduced in  Section~\ref{sec:inival}
    are constructed from point delays: From
    $\lambda=s\sqrt{a}-\alpha$, $\alpha\in\Complexes$ we obtain
    $2C(x)=\mathrm{e}^{s\sqrt{a}-\alpha}+\mathrm{e}^{-(s\sqrt{a}-\alpha)}$. Note that
    in this case our results are simply a restatement of those presented in
    \cite{BreLoi96msca,Gluesing97siam}. 
  \end{remark}

\subsubsection{The polynomial $\lambda^2-\sigma$ is irreducible over $k$}\label{sec:irreducible}
As indicated by Example \ref{ex:ideals}, the second case, \ie, the
equation $\lambda^2-\sigma=0$ has no solutions over $k$, is much more
challenging than the first one. There, the ring $\RingG{\Naturals}$
corresponds basically to the ring $\Rationals[x,y]/[x^2+y^2-1]$ of
trigonometric polynomials which is obtained from $\RingG{\Naturals}$
for $\sigma=-1$ and $k=\Rationals$. The latter ring is lacking the
pleasing properties of a PID or even a \Bezout\
domain\footnote{Actually, the trigonometric ring is a Dedekind domain
  \cite{Cohn2003ba,Cohn2003fa,Picavet03trigon}.}.  However, 
Example \ref{ex:ideals} suggests, that the difficulties can be
circumvented when allowing to halve the argument, \ie, working with
$\RingG{\Rationals}$ instead of $\RingG{\Naturals}$.

\begin{definition}
  For any nonzero $r\in \RingG{\domain}$ the \emph{norm} $\nu(r)$ is defined
  as the highest $\alpha\in \domain^+$ such that at least one of the
  coefficients $a_{\alpha}$ and $b_{\alpha}$ in \eqref{eq:poly} is
  nonzero.
\end{definition}
\begin{lemma}\label{lmm:division_step}
  Let $S$ the multiplicative subset of $\RingG{\Naturals}$ consisting
  of all the elements such that either the coefficients with odd or
  those with even indices vanish. More precisely, any element $s$ of $S$ can be written as
  \begin{align*}
    s=\sum_{i\in I_s} a_{s,i} C_{i}+b_{s,i} S_{i}, \quad I_s=\left\{\nu(s)-2i \,\big|\, i\in\Integers,\,0\leqslant i\leqslant \frac{\nu(s)}{2} \right\}
  \end{align*}
  Let $p,q\in S$ the norms of which are strictly positive. Without
  loss of generality assume $\nu(p)\geqslant \nu(q)$. Consider the ideal
  $\mathfrak{I}=(p,q)$ generated by $p$ and $q$.  Then there exists
  $\bar p,\bar q\in S$ with $\mathfrak{I}=(\bar p, \bar q)$ and either
  $\nu(p)>\nu(\bar p)\geqslant \nu(\bar q)$ or $\bar q=0$.
\end{lemma}
\begin{proof}
  In the following, three different cases are considered.

  \textbf{Case 1.} If $\nu(p)>\nu(q)$ one can apply a division step
    similar to that of polynomials. More precisely, we will show that
    there exists $r,h\in S$ with either $r=0$ or $\nu(r)<\nu( p)$ such that $p=qh+r$.
    Then we may set $\bar p=q$, $\bar q= r$ (or vice versa) to complete the discussion of the first case.
    
    In order to show that $r$, $h$ with the claimed properties exist set 
   $$
    h=a_{h} C_\Delta+b_{h} S_\Delta,\quad\Delta={\nu(p)-\nu(q)}
   $$ 
   where the coefficients $a_h,b_h \in \field$ have to be determined appropriately.
   It follows 
  \begin{align*}
      s=hq&=\sum_{i\in I_q}\Big((a_{h} a_{q,i} C_{i}C_{\Delta}+b_{h} a_{q,i} C_{i}S_{\Delta})+(a_{h} b_{q,i} S_{i}C_{\Delta}+b_{h} b_{q,i} S_{i} S_{\Delta})\Big)\\
&=\frac{1}{2\sigma}\sum_{i\in I_q} \Big((\sigma a_{h} a_{q,i}+b_{h} b_{q,i}) C_{\Delta+i}+(\sigma a_{h} a_{q,i}-b_{h} b_{q,i})C_{\Delta-i}\\
             &\qquad\qquad+\sigma(b_{h} a_{q,i}+ a_{h} b_{q,i}) S_{\Delta+i}+ \sigma(b_{h} a_{q,i}-a_{h} b_{q,i}) S_{\Delta-i}\Big)\\
&=\sum_{i\in I_p}   a_{s,i} C_{i}+b_{s,i} S_{i}
\end{align*}
where  the leading coefficients are given by
$$
a_{s,\nu(p)}=\frac{1}{2\sigma}\big(\sigma a_{h} a_{q,\nu(q)}+b_{h} b_{q,\nu(q)}\big),\quad
b_{s,\nu(p)}=\frac{1}{2}\big(b_{h} a_{q,\nu(q)}+ a_{h} b_{q,\nu(q)}\big).
$$
From this equation and from $r=hq-p$ the norm of $r$ is smaller than
that of $p$ if and only if $a_{h}$, $b_h$ satisfy
\begin{equation}\label{eq:lmm:glsys}
    \begin{pmatrix}
      a_{q,\nu(q)}&\sigma^{-1}{b_{q,\nu(q)}}\\
      b_{q,\nu(q)}&a_{q,\nu(q)}
    \end{pmatrix}
    \begin{pmatrix}
     a_{h}\\
     b_{h} 
    \end{pmatrix}=
    2\begin{pmatrix}
      a_{p,\nu(p)}\\
      b_{p,\nu(p)}
    \end{pmatrix}.
  \end{equation}
By the definition of the norm  at least one of the coefficients $a_{q,\nu(q)}$ and $a_{q,\nu(q)}$ is nonzero. Since, additionally, $\sqrt{\sigma}\not\in \field$ it follows
$\sigma a_{q,\nu(q)}^2- b_{q,\nu(q)}^2\ne 0$ and  $a_h, b_h$ can be alway chosen according to \eqref{eq:lmm:glsys}.

\textbf{Case 2.}  If $\nu(p)=\nu(q)$ and for some
$c\in \field$ the equations $a_{q,\nu(q)}=c a_{p,\nu(p)}$, $b_{q,\nu(q)}=c b_{p,\nu(p)}$ hold,
the ideal $\mathfrak{ I}$ is generated by $\bar p=p$, $\bar q=q-cp$ where
$\nu(\bar q)<\nu(\bar p)$. If $\bar q=0$ the proof is complete otherwise we
can proceed according to the first case with the pair $\bar
p,\bar q$ instead of $p,q$.

\textbf{Case 3.} If $\nu(p)=\nu(q)$ but we are not in the  second case set 
  \begin{subequations}
    \begin{align}
      \begin{pmatrix}
        p&q
      \end{pmatrix}^T
      &=A_1
      \begin{pmatrix}
        \tilde p&  \tilde q
      \end{pmatrix}^T\label{eq:A1}\\
      \begin{pmatrix}
        \tilde p&\tilde q
      \end{pmatrix}^T
      &=A_2
      \begin{pmatrix}
        \bar p& \bar q
      \end{pmatrix}^T\label{eq:A2}
    \end{align}
  \end{subequations}
with
  \begin{align*}
    A_1=\begin{pmatrix}
      a_{p,n}&b_{p,n}\\
      a_{q,n}&b_{q,n}
    \end{pmatrix},\quad A_2=
    \begin{pmatrix}
      C_1& \sigma S_1\\
       S_1& C_1
    \end{pmatrix},\quad n=\nu(q)=\nu(p)
  \end{align*}
  Obviously, $p,q$ belong to the ideal generated by $\bar p,\bar q$.
  Both matrices, $A_1$ and $A_2$, are invertible, the first one since otherwise
  we would be in the second case, the latter one since, by
  (\ref{eq:hypfun:propA}), its
  determinant equals $1$.  Thus, $(\bar p, \bar q)=(\tilde p,\tilde
  q)=(p,q)$.

  It remains to show that the norms of $\bar p$ and $\bar q$ are
  both smaller than $n$. From equation (\ref{eq:A1}) one obtains
  $\nu(\tilde p)=\nu(\tilde q)=n$ with $a_{\tilde p,\nu(q)}=b_{\tilde
    q,\nu(q)}=1$, $b_{\tilde p,\nu(q)}=a_{\tilde q,\nu(q)}=0$.
  From (\ref{eq:A2}) one has
    \begin{align*}
      \bar p&=C_1C_{n}-\sigma S_1S_{n}+\sum_{i\in I_p^*}
      C_1\big(a_{\tilde p,i}C_i+b_{\tilde p,i}S_i\big)-
      \sigma S_1\big(a_{\tilde q,i}C_i+b_{\tilde q,i}S_i\big)\\
      \bar q&=C_1S_{n}-S_1C_{n}+
      \sum_{i\in I_p^*}
      C_1\big(a_{\tilde q,i}C_i+b_{\tilde q,i}S_i\big)-
      S_1\big(a_{\tilde p,i}C_i+b_{\tilde p,i}S_i\big)
    \end{align*}
with $I_p^*=I_p\setminus\{n\}$.
The  norms of the sums in the above expression are at most $n-1$
while for the leading terms one obtains according to
 (\ref{eq:hypfun:prop:1}) 
$$
C_1C_n-\sigma S_1S_n=C_{n-1},\quad C_1S_n-S_1C_n=S_{n-1}.
$$
Thus, the norms of $\bar p$, $\bar q$ cannot exceed $n-1$.
\end{proof}

\begin{lemma}\label{lemma:division_algorithm}
  Let $p,q\in S\subset \RingG{\Naturals}$ with $\nu( p)\geqslant \nu( q)$. Then there exists $\bar p, \bar q\in S\cap (p,q)$ such that
  $(p,q)=(\bar p,\bar q)$ and $\nu(\bar q)<\nu(q)$, $\nu(p)\leqslant\nu(q)$ or $\bar q=0$.
\end{lemma}
\begin{proof}
  By Lemma \ref{lmm:division_step} $(p,q)=( p^*, q^*)$ with
  $\nu(p)>\nu(p^*)\geqslant \nu(q^*)$ or $q^*=0$. In the latter case the
  claim has been proved.  Otherwise, repeat the above argument $p^*,
  q^*$ until we are in the claimed situation which happens after at
  most $\nu(p)-\nu(q)+1$ steps.
\end{proof}

\begin{proposition}\label{prop:principal_ideals_in_Rn}
Any ideal $\mathfrak{I}$ in $\RingG{\Naturals}$ generated by a subset $\mathfrak{G}$ of $S$ is principal.
\end{proposition}
\begin{proof}
  \textbf{Step 1.} We show that up to multiplication with units there
  is only one element $q$ of lowest norm $\nu(q)=n$ in $S\cap
  \mathfrak{I}$.  To this end, assume there are at least two such
  elements, say $p$ and $q$.  By Lemma \ref{lemma:division_algorithm}
  there exist $\bar p,\bar q\in S$ with $(\bar p,\bar q)=(p,q)$ where
  $n>\nu(\bar p)\geqslant\nu(\bar q)$ or $n\geqslant\nu(\bar p)$ and
  $\bar q=0$.  Since $n$ is the lowest possible norm for an element of
  $\mathfrak{I}\cap S$, only the case $n=\nu(\bar p)$ and $\bar q=0$
  remains. But this can happen only if we are in case~2 of
  Lemma~\ref{lmm:division_step} having $\bar p=p$ and $q=cp$, $c\in
  \field^{\times}$.
  
  \textbf{Step 2.} We now show that any element of $\mathfrak{G}$
  belongs to $(q)$ where $q$ is defined as in the first step. To this
  end chose any element $p$ from $\mathfrak{G}$. Applying case $1$ of
  Lemma~~\ref{lmm:division_step} several times one gets $p=hq+r$,
  $\nu(r)\leqslant n$, $r\in S$. Since, by assumption, $q$ has the smallest
  possible norm, it follows $\nu(r)=n$ or $r=0$. This in turn yields
  $r=cq$, $c\in \field$ according to Step~1. Finally, we have $p=(h+c)q$ and, therefore,
  $\mathfrak{I}=(q)$.
\end{proof}
\begin{proposition}\label{prop:bezout2}
Any finitely generated ideal  $\mathfrak{I}$ in $\RingG{\Rationals}$ is principal, \ie, $\RingG{\Rationals}$ is 
a B\'ezout domain.
\end{proposition}
\begin{proof}
  Let  $\mathfrak{I}=(r_1,\dots,r_m)$ for some $m\in\Naturals$. Write the generators according  to (\ref{eq:poly}), \ie,
   \begin{equation}\label{eq:poly:1}
     r_j=\sum_{i=0}^{n_j} a_{\alpha_{j,i}} C_{\alpha_{j,i}} +     b_{\alpha_{j,i}} S_{\alpha_{j,i}},\quad n_j\in\Naturals, \quad a_{\alpha_{i,j}},b_{\alpha_{i,j}}\in \field, \quad \alpha_{i,j}\in Q^+.
   \end{equation}
   Let $d$ be a common denominator of all the $\alpha_{i,j}$ occurring in these
   equations. Then the generators of $\mathfrak{I}$ can be identified with
   elements of the subset $S$ defined in Lemma~\eqref{lmm:division_step} of the
   Ring $R_{\Naturals}$ via the embedding $E:R_{\Naturals}\rightarrow
   R_{\Rationals}$ which is defined by $\tilde C_{2}\mapsto C_{1/d}$, $\tilde
   S_{2}\mapsto S_{1/d}$. Let $\tilde r_1,\dots,\tilde r_m$ elements of $R_N$
   the images of which are $r_1,\dots,r_m$. The ideal $\tilde{\mathfrak{I}}$
   generated by $\tilde r_1,\dots,\tilde r_m$ is principal by
   Proposition~\ref{prop:principal_ideals_in_Rn}. Consequently, $\mathfrak{I}$
   is generated by the image of the generator of $\tilde{\mathfrak{I}}$ under $E$.
\end{proof}

\begin{remark}\label{rem:nopid}
  Note that neither $\RingG{\Rationals}$ nor $\RingG{\Naturals}$ are 
  principal ideal domains (PID).  The first is not Noetherian: As an
  example for an ideal that is not finitely generated take
  $(\{S_{1/2^n}|n\in\Naturals\})$.  Moreover, $\RingG{\Naturals}$ is not a
  PID since there are finitely generated ideals that cannot be
  generated by one single element: The ideal $(S_1,C_1+1)$ viewed as
  an element of $\RingG{\Rationals}$ is generated by $C_{1/2}$ which does
  not belong to $\RingG{\Naturals}$. 
\end{remark}

\subsection{$\Ring_\Rationals$ is a \Bezout\ domain}\label{sec:bezout}
We are now in position to prove that $\Ring_\Rationals$ is a \Bezout\ domain.
After the preparation done in the previous subsection the remaining steps are
very similar to those given in \cite{BreLoi96msca,Gluesing97siam}. In
particular, the proof of Lemma~\ref{lemma:coprime} which prepares 
Theorem~\ref{theorem:bezout} is strongly inspired by
\cite[Theorem~1]{BreLoi96msca}.

\begin{lemma}\label{lemma:coprime}
  For two coprime elements $p,q\in\Ring_{\Rationals}$ there exist $a,b\in\Ring_{\Rationals}$ such that
  $ap+bq=1$.
\end{lemma} 
\begin{proof}
  By Prop.~\ref{prop:bezout2} (resp.\ Cor.~\ref{corollary:bezout1}) $\RingS{\Rationals}$ is
  a \Bezout\ domain.  Thus, there exist $a,b\in\Complexes[s,\mathfrak{S}]\subset\Ring_\Rationals$  such
  that $ap+bq=h$ where $h\in \Complexes[s]$.  Write $h$ as product
  $h=\prod_{i=1}^N(s-s_i)$ and proceed by induction (we do not assume $s_i\ne
  s_j$).

  Assume there exist $a,b\in\Ring_{\Rationals}$ with
  $ap+bq=\prod_{i=1}^N(s-s_i)$.  In the following, for any
  $\gamma\in\MikComp$ we set $\bar\gamma=\LaplaceTrans(\gamma)(s_N)$,
  with the entire function $\LaplaceTrans(\gamma)$ being the Laplace
  transform of $\gamma$ given according to
  App.~\ref{sec:miku_divis}. By the coprimeness of $p,q$ and by
  Lemma~\ref{thrm:mik:teilbarkeit} in the Appendix, $a^\star$ and
  $b^\star$, defined by
  \begin{equation*}
    \begin{split}
          a^\star&=\begin{cases}
                  \dfrac{\bar q a-q \bar a }{\bar q(s-s_N)},&  (s-s_N)\nmid    q\\
                  \dfrac{a}{s-s_N},&  (s-s_N)\mid q
                 \end{cases},\quad
          b^\star=\begin{cases}
                  \dfrac{\bar p b-p \bar b }{\bar p(s-s_N)},&  (s-s_N)\nmid    p\\
                  \dfrac{b}{s-s_N},&  (s-s_N)\mid p
                 \end{cases}
    \end{split}
  \end{equation*}
  belong to $\Ring_\Rationals$. One easily verifies that $pa^\star+qb^\star=\prod_{i=1}^{N-1}(s-s_i)$.
  Applying this step $N$ times completes the proof.
\end{proof}
\begin{theorem}\label{theorem:bezout}
  The ring $\Ring_{\Rationals}$ is a B\'ezout domain, \ie, any finitely generated
  ideal is principal.
\end{theorem} 
\begin{proof}
  We show that any two elements $p,q\in \Ring_\Rationals$ possess a
  common divisor $\tilde c$ that can be written as linear
  combination of $p$, $q$. (It is then the unique greatest common divisor of $p$ and $q$.)  

According to section
  \ref{sec:ring:field} the ring $\RingS{\Rationals}$ is a \Bezout\
  domain.  Consequently, there are elements $a,b\in
  \Complexes[s,\mathfrak{S}_{\Rationals}]$ such that
  \begin{equation}
    \label{eq:h1}
       c=ap+bq\in  \Complexes[s,\mathfrak{S}_{\Rationals}]
  \end{equation}
  is a g.c.d.\ in
  $\RingS{\Rationals}$.  Hence, $p/c$ and $q/c$ belong to
  $\RingS{\Rationals}$. In particular, there are $n_i\in\Ring_{\Rationals}$,
  $d_i\in\Complexes[s]$ with $\gcd_{\Ring_{\Rationals}}(n_i,d_i)=1$
  ($i=1,2$) such that $p/c=n_1/d_1$ and $q/c=n_2/d_2$. It follows
  $p d_1=c n_1$,  $q d_2=c n_2$.
  Consequently, both $d_1$ and $d_2$ divide $c$ in $\Ring_{\Rationals}$.
  Since $d_1$ and $d_2$ are polynomials, they possess a least common
  multiple $h=d_1d_2/\gcd(d_1,d_2)$, and it follows $\tilde c=c/h\in \Ring_{\Rationals}$.
  Clearly, $\tilde c$ divides both, $p$ and $q$.  
  Dividing (\ref{eq:h1}) by $\tilde c$ yields the equation
  \begin{equation}
  a \underbrace{n_1 d_2/\gcd(d_1,d_2)}_{\tilde p}+b \underbrace{n_2 d_1/\gcd(d_1,d_2)}_{\tilde q}=\underbrace{d_1d_2/\gcd(d_1,d_2)}_{h}.\label{eq:h2}
\end{equation}
  By the coprimeness of $n_1$ and $d_1$, resp.\ $n_2$ and $d_2$, it follows $\gcd(\tilde p,h)=d_2/\gcd(d_1,d_2)$, resp.\
  $\gcd(\tilde q,h)=d_1/\gcd(d_1,d_2)$. Thus, $\gcd(\tilde q,h)$ and $\gcd(\tilde p,h)$ are coprime
  and, since by equation \eqref{eq:h2} any common divisor of $\tilde p$ and $\tilde q$ divides $h$, we can finally conclude the coprimeness of $\tilde p$ and $\tilde q$.
  Thus, by Lemma~\ref{lemma:coprime}, there are $a^\star, b^\star\in\Ring_\Rationals$ such that
   $a^\star \tilde p+b^\star \tilde q=1$. The claim follows directly by multiplying this equation 
   by $\tilde c$.
\end{proof}
\newpage
\section{Controllability analysis}\label{sec:controllability}
\subsection{Systems and Dynamics}
\begin{notation}
  Unless otherwise stated, the submodule spanned by a subset $S$ of an
  $R$-module $M$ is written $[S]$.
\end{notation}

\begin{definition}
  An $R$-\emph{system} $\Lambda$, or a \emph{system over} $R$, is an
  $R$-module.
\end{definition}
\begin{definition}
  A \emph{presentation matrix} of a finitely presented $R$-system $\Sigma$ is a matrix $P$ such that $\Sigma\cong [v]/[Pv]$ where
  $[v]$ is free with basis $v$.
\end{definition}

\begin{definition}
  An $R$-\emph{dynamics}, or a \emph{dynamics over} $R$, is an
  $R$-system $\Lambda$ equipped with an \emph{input}, \ie, a subset
  $\Vect{u}$ of $\Lambda$ which may be empty, such that the quotient
  $R$-module $\Lambda/[\Vect{u}]$ is torsion.  The input $\Vect{u}$ is
  \emph{independent} if the $R$-module $[\Vect{u}]$ is free, with
  basis $\Vect{u}$.
\end{definition}

\begin{definition}
  An \emph{output} $\Vect{y}$ is a subset, which may be empty, of
  $\Lambda$.  An \emph{input-output $R$-system}, or an
  \emph{input-output system over} $R$, is an $R$-dynamics equipped
  with an output.
\end{definition}

\begin{definition}
  Let $A$ be an $R$-algebra and $\Lambda$ be an $R$-system.  The
  $A$-module $A \otimes_R \Lambda$ is an $A$-system, which
  \emph{extends} $\Lambda$.
\end{definition}

\subsection{System controllabilities}\label{subsecAContr}
In this section we emphasize several controllability notions which are defined
directly on the basis of the above system definition without referring to a
solution space. For the latter we refer to the next subsection. Let us start
with some purely algebraic definitions:
\begin{definition}[\emph{see} \cite{FliMou98cocv}]
  Let $A$ be an $R$ algebra. An $R$-system $\Lambda$ is said to be
  \emph{$A$-torsion free controllable} (resp.~\emph{$A$-projective
    controllable}, \emph{$A$-free controllable}) if the $A$-module
  $A\otimes_R\Lambda$ is torsion free (resp.~projective, free). An
  \emph{$R$-torsion free} (resp.~\emph{$R$-projective},
  \emph{$R$-free}) \emph{controllable} $R$-system is simply called
  \emph{torsion free} (resp.~\emph{projective}, \emph{free})
  \emph{controllable}.
\end{definition}
Elementary homological algebra (\seeeg{} \cite{Rotman1979}) yields
\begin{proposition}
  $A$-free (resp.~$A$-projective) controllability implies
  $A$-pro\-jec\-tive (resp.~$A$-torsion free) controllability.
\end{proposition}
\begin{proposition}
  $R$-free  controllability implies
  $A$-free  controllability for any $R$-algebra $A$.
  More generally, given any $R$-system $\Sigma$ that is a direct sum of a torsion module $\mathrm{t}\Sigma$ and a free module
  $\Lambda$, the extended system $A\otimes_R\Sigma$ is a direct sum of 
  the torsion module $A\otimes_R\mathrm{t}\Sigma$ and the free module $A\otimes_R\Lambda$.
\end{proposition}

The importance of the notions of torsion free and free controllability
is intuitively clear: While the first one refers to the absence of a
nontrivial subsystem which is governed by an autonomous system of
equations, the latter refers to the possibility to freely express all
system variables in terms of a basis of the system module. For this
reason, and, secondarily, in reminiscence to the theory of nonlinear
finite dimensional systems, we have the following:
\begin{definition}
  Take an $A$-free controllable $R$-system $\Lambda$ with a finite output
  $\Vect{y}$. This output is said to be $A$-\emph{flat}, or $A$-\emph{basic},
  if $\Vect{y}$ is a basis of $A\otimes_R\Lambda$. If $A \cong R$ then 
  $\Vect{y}$ is simply called \emph{flat}, or \emph{basic}.
\end{definition}

In finite dimensional linear systems theory, the so called Hautus
criterion is a quite popular tool for checking controllability. This
criterion has been generalized to delay systems (\seeeg
\cite{Mounier1998fm}) and to the more general convolutional systems
defined over $\DistComp$ \cite{VetZam00mtns} and $\MikComp$
\cite{Woi07diss}. All those rings may be embedded into the ring of
entire functions via the Laplace transform. This motivates the
following quite general definition:
\begin{definitionandproposition}\label{def:spectr}
    Let $R$ be any ring that is isomorphic to a subring of the ring $\entire$ of
    entire functions with pointwise defined multiplication. Denote the embedding
    $R\rightarrow\entire$ by $\LaplaceTrans$.  A finitely presented $R$-system
     with presentation matrix $P$ is said to be spectrally
    controllable if one of the following equivalent conditions holds:
    \begin{enumerate}\renewcommand{\theenumi}{(\roman{enumi})}
    \item The
      $\entire$-matrix $\hat P=\LaplaceTrans(P)$ satisfies $\exists
      k\in\Naturals:\forall\sigma\in\Complexes:\rank\hat P(\sigma)=k$.
    \item The module $\Sigma_\entire=\entire\otimes_R \Sigma$ is free.
    \end{enumerate}
  \end{definitionandproposition}
  \begin{proof}
    The ring $\entire$ is an elementary divisor domain
    \cite{Helmer1943}. As a consequence, over $\entire$, any matrix
    admits a Smith normal form by left and right multiplication with
    unimodular matrices. Since the units in $\mathcal{O}$ are just
    the functions which possess no zeros in $\Complexes$, the rank
    of the Smith normal form equals the rank of $\hat P(\sigma)$ for
    any $\sigma\in\Complexes$. Thus, the rank of $\hat P(\sigma)$ remains
    constant if and only if the nonzero entries of the Smith normal
    form possess no zeros in $\Complexes$ which, in turn, is
    equivalent to the absence of nontrivial torsion elements in
    $[v]/[\hat Pv]$, \ie, in $\Sigma_{\entire}$.
  \end{proof}
  \begin{proposition}\label{prop:spectr}
    Let $R$ be any \Bezout\ domain that is isomorphic to a subring of
    $\entire$ with the embedding $R\rightarrow\entire$ denoted by
    $\LaplaceTrans$.  Then the notions of spectral controllability and
    $R$-torsion free controllability are equivalent if and only if
    $\LaplaceTrans$ maps non-units in $R$ to non-units in $\entire$.
  \end{proposition}
  \begin{proof}
    Since $R$ is a \Bezout\ domain, torsion freeness of $\Sigma$
    implies freeness. Tensoring with the free module $\entire$ yields
    another free module $\Sigma_\entire$, and, thus, by Definition and
    Proposition \ref{def:spectr}, spectral controllability.  Again,
    since $R$ is a \Bezout\ domain, any presentation matrix admits a
    Hermite form. Thus, the torsion submodule $\mathrm{t}\Sigma$ of
    $\Sigma$ can be presented by a triangular square matrix
    $\mathrm{t}P$ of full rank.  If $\Sigma$ is not torsion-free, at
    least one diagonal entry of this matrix is not a unit in $R$. If
    this entry is mapped to a non-unit in $\Sigma_\entire$ by
    $\LaplaceTrans$, it admits a complex zero $\sigma_0$. Thus,
    $\LaplaceTrans(\mathrm{t}P)$ has a loss off rank at
    $\sigma=\sigma_0$. Contrary, if there is a non unit $r\in R$ which
    corresponds to a unit $\hat r\in\entire$, consider
    $\Sigma\cong[\tau]/[r\tau]$.  Obviously, the image of $\tau$ in
    $\Sigma_\entire$ is zero. Thus, the trivial module
    $\Sigma_\entire$ is torsion free.
  \end{proof}
  \begin{remark}
    Note that, under the additional assumption that $\Sigma$ admits a
    presentation matrix of full row-rank, the assumption of $R$ being a \Bezout\
    domain may be replaced by a less restrictive one. In this case, equivalence
    of $(\mathcal{Q}\otimes_R R)\cap\entire$-torsion free controllability, with
    $\mathcal{Q}$ the ring of rational functions in one complex variable, and
    spectral controllabillity may be established (\seeeg\ \cite{Mounier1998fm,Woi07diss} for different examples).
  \end{remark}
  We are now able to state the main result of our paper:
\begin{theorem}\label{thrm:main}
  The convolutional system $\Sigma$ defined in
  Definition~\ref{def:Sigma} is free if and only if it is torsion
  free. More generally
  $\Sigma=\mathrm{t}\Sigma\oplus\Sigma/\mathrm{t}\Sigma$ where
  $\mathrm{t}\Sigma$ is torsion and $\Sigma/\mathrm{t}\Sigma$ is
  free. Moreover, $\Sigma$ is spectrally controllable if and only if
  it is torsion free.
\end{theorem}
\begin{proof}
  Recall that, according to Definition~\ref{def:Sigma},
  \mbox{$\Sigma\cong\Ring\otimes_{\Ring_\Rationals}\!\Sigma_\Rationals$} and
  $\Ring_\Rationals$ is a \Bezout\ domain by Proposition~\ref{theorem:bezout}. Since
  the first assertion holds for finitely presented modules over any \Bezout\
  domain, it holds for $\Sigma_\Rationals$.  The second
  assertion follows from Proposition~\ref{prop:spectr}. (The fact 
  that the Laplace transform maps any non-unit of $\Ring_{\Rationals}$ to a
  non-unit in $\entire$ is obvious.) Clearly, both results hold as
  well for $\Sigma$, which is obtained by an extension of scalars.
\end{proof}

\subsection{Trajectorian controllability}\label{sec:traj-contr}
In this section we will give two different interpretations of our
algebraic controllability results that directly refer to trajectories of the
system, \ie, to (generalized) functions which may be assigned to the
system variables. To this end we need to introduce the notions of a \emph{solution
space} and a \emph{trajectory}.
\begin{definition}\label{def:solspace}
  Let $\Sigma$ be an $R$-system and $\mathscr{F}$ a space of generalized functions.
  The space $\mathscr{F}$ is called  
  a \emph{solution space} of $\Sigma$ if 
  it can be equipped with the structure of an $R$-module.
\end{definition}
\begin{definition}[see \cite{FliMou02ima}]\label{def:trj}
  Let $\mathscr{F}$ be a solution space of an $R$-system $\Sigma$. An
  $\mathscr{F}$-\emph{trajectory} of $\Sigma$ is an
  element of $\Hom_R(\Sigma,\mathscr{F})$.
\end{definition}

The crux of the first controllability notion (Def.~\ref{def:TrajContr}) is the possibility to
assign an arbitrary (generalized) function from $\mathscr{F}$ to \emph{any} system variable.
\begin{definition}[see \cite{FliMou02ima}]
\label{def:TrajContr}
  An  $R$-system is called $\mathscr{F}$-\emph{trajectory control\-lable} if for any element
  $a\in\Sigma$ and any $b\in\mathscr{F}$ there exists a trajectory $f$ with $f(a)=b$.
\end{definition}

The following result is borrowed from \cite{FliMou02ima} and applies
to any torsion-free controllable $R$-system where $R\subset \Mik$.
\begin{proposition}
  The system $\Sigma_\Reals/\mathrm{t}\Sigma_\Reals$ is $\Mik$-trajectory-controllable. 
\end{proposition}

Another controllability notion is the following due to \cite{Wil91}. As the 
above it relies on the notion of a trajectory. However, since it refers to the
possibility of connecting trajectories, the notions of future and past come into
play. Thus, an appropriate solution space should allow the definition of such
local properties. This is not possible for the field of \Mikusinski\ operators
in general but for its subring $\MikReg$ and more generally for the 
 space $\Boem$ of Boehmians
\cite{Boehme1973ams,Mikusinski1990amh}. The controllability criterion in the
behavioural framework is the possibility of concatenating trajectories.  In our
algebraic setting we may formulate this criterion as follows.
\begin{definition}[cf.\ \cite{Wil91,RocWil97}]
\label{def:BehavContr}
  Let $\Sigma$ be an $R$-system and $\mathscr{F}$  a solution space of $\Sigma$ that
  possesses the structure of a sheaf on $\Reals$.  Then $\Sigma$ is called
  \emph{$\mathscr{F}$-behavioral-controllable} if for any two trajectories
  $f_1,f_2\in\Hom(\Sigma,\mathscr{F})$ there exists $f\in
  \Hom(\Sigma,\mathscr{F})$ such that for any $a\in\Sigma$ there are
  $t_1^a,t_2^a\in\Reals$ with
  $f(a)|_{(-\infty,t_1^a)}=f_1(a)|_{(-\infty,t_1^a)}$ and
  $f(a)|_{(t_2^a,\infty)}=f_2(a)|_{(t_2^a,\infty)}$.
\end{definition}
\begin{theorem}
  The system $\Sigma_\Reals/\mathrm{t}\Sigma_\Reals$, where $\Sigma_\Reals$ is defined in Definition \ref{def:Sigma}, is
  $\Boem$-behavioural controllable.
\end{theorem}
\begin{proof}
  Since $\Sigma_\Reals/\mathrm{t}\Sigma_\Reals$ is free, any homomorphism is uniquely
  determined by the functions assigned to the basis.  Thus, for the basis
  $b=b_1,\dots,b_n$  we may chose $t_1^b>t_2^b$ and set
  $$
    f(b)=\begin{cases}
           f_1(b),&t<t_1^b\\
           f_2(b),&t>t_2^b
         \end{cases}
  $$
  Moreover, any $a\in\Sigma_\Reals/\mathrm{t}\Sigma_\Reals$ is given by
  $a=\sum_{i=0}^n\alpha_ib_i$ where the $\alpha_i$ have compact support.  Thus there
  exist $T_1$, $T_2$ such that $\supp\alpha_i \subseteq [T_1,T_2]$, $i=1,\dots,n$.
  The claim follows by an application of the theorem of supports 
   $
       t^a_1=t_1^b+T_1, \quad t^a_2=t_2^b+T_2$ (see.\ \cite{Boehme1973ams}).
\end{proof}
\begin{remark}
  When distinguishing the cases $a>0$ and $a=0$ in (\ref{eq:chareq})
  one could alternatively prove $\Cinfty$-behavioural controllability
  (resp.\ $\Dist$-behavioural controllability) in the case $a>0$ or
  $\Gevrey{2}$-behavioural controllability (resp.\
  $\DistGevrey{2}$-behavioural controllability) in the case $a=0$,
  where $\Cinfty$ is the space of infinitely differentiable
  functions, $\Dist$ the space of Schwartz-Distributions,
  $\Gevrey{2}$ the space of Gevrey-Functions of order less than
  $2$, and $\DistGevrey{2}$ the space of Gevrey ultradistributions.
\end{remark}
\section{An example: two boundary coupled equations}\label{sec:example}
In order to illustrate our results, in the following we discuss a
simple example. Consider the system of two second order
equations 
\begin{subequations}\label{EqEdpSecondOrdre}
\begin{equation}\label{EqEdpSecondOrdre:pde}
  \partial_x^2 w_i(x)   = \sigma w,\quad i=1,2,
\end{equation}
defined on an open neighbourhood $\Omega_i$ of
$[0,\ell_i]\subset\Reals$,
where $\sigma=\alpha s^2+\beta s+c$. Those equations are coupled via
the boundary conditions ($i=1,2$)
\begin{align}\label{EqEdpSecondOrdre:bc}
\mu_{i1} w_i(\ell_i)+ \mu_{i2} w_i'(\ell_i) &= 0 \\
w_i(0)&=u.
\end{align}
\end{subequations}

According to Section~\ref{sec:inival}, the general solution of the
initial value problems associated with \eqref{EqEdpSecondOrdre:pde}
reads ($i=1,2$)
\begin{align}
  \label{eq:networks:GenSol}
  \begin{pmatrix}
  w_i(x)\\w_i'(x)
  \end{pmatrix}=
  \begin{pmatrix}
    C(x-\ell_i)&S(x-\ell_i) \\
    \sigma S(x-\ell_i)& C(x-\ell_i)
  \end{pmatrix}
  \begin{pmatrix}
    \varK_{i1}\\\varK_{i2}
  \end{pmatrix},
\end{align}
with $\varK_{i1}=w_i(\ell_i)$,
$\varK_{i2}=\partial_x w_i({\ell_i})$.
The boundary conditions at $x=\ell_i$ yield
\begin{subequations}\label{eq:ex:presentation1}
  \begin{align}
    \label{eq:ex:presentation1:1}
    \mu_{i1}\varK_{i1}+\mu_{i2}\varK_{i2}&=0\\
    \label{eq:ex:presentation1:2}
    C(\ell_i)\varK_{i1}-S(\ell_i)\varK_{i2}&=u.
  \end{align}
\end{subequations}
Here, the relations $S(-\ell_i)=-S(\ell_i)$ and $C(-\ell_i)=C(\ell_i)$,
derived in Section \ref{sec:inival}, have already been incorporated.

Thus, according to Definition \ref{def:Sigma}, the convolutional
system $\Sigma$ associated with the boundary value problem
(\ref{EqEdpSecondOrdre}) is the $\Ring$ module
$[\varK_{11},\varK_{12},\varK_{21},\varK_{22},u]$ the generators of which are
subject to the equations (\ref{eq:ex:presentation1}).

In order to reduce the number of equations, we aim to introduce new
variables $\omega_1$ and $\omega_2$ such that
(\ref{eq:ex:presentation1:1}) is satisfied automatically, \ie,
\begin{equation*}
  \varK_{i1}=-\mu_{i2}\omega_i,\quad
  \varK_{i2}= \mu_{i1}\omega_i,\quad i=1,2.
\end{equation*}
Indeed, since
\begin{equation*}
\omega_i=\dfrac{1}{\mu_{i1}^2 + \mu_{i2}^2}\left(-\mu_{i2}\varK_1+\mu_{i1}\varK_2\right),\quad i=1,2,
\end{equation*}
the new variables belong to $\Sigma$. Using the new generators $\omega_1$, $\omega_2$ and $u$,
equation (\ref{eq:ex:presentation1:2}) may be rewritten to obtain
\begin{equation*}
   u=-p_i\omega_i,\quad p_i=\mu_{i2}C(\ell_i)+\mu_{i1}S(\ell_i),\quad i=1,2.
\end{equation*}
Thus, $p_1\omega_1-p_2\omega_2=0$, and $\Sigma\cong[\tilde\omega_1,\tilde\omega_2]/[p_1\tilde\omega_1-p_2\tilde\omega_2]$.

In accordance with  Section \ref{sec:model} assume
that $\ell_i=n_i\ell$, with $n_i\in\Naturals$ and $i=1,2$.  Thus, by
Theorem~\ref{thrm:main}, checking spectral, torsion free, and free
controllability are equivalent. Since the aim of this section is not
the presentation of a general controllability analysis for the boundary value
problem 
(\ref{EqEdpSecondOrdre}) but rather to give an example for the
application of the derived algebraic results, we shall restrict
ourselves to particular values for $n_1$ and $n_2$. In order to avoid
tedious computations, we chose simply $n_1=1$, $n_2=2$. Apart from that,
we  discuss the generic case only, \ie, we do not care about singularities
which may occur for particular values of the $\mu_{ij}$, $i,j=1,2$.

Applying the algorithms of Section \ref{sec:ring:field} we obtain
$p_1r_1+p_2 r_2=\epsilon$ with
\begin{align*}
  r_1&=2\big((\mu_{21}\mu_{11}-\mu_{22}\mu_{12}\sigma)C(\ell)+\left(\mu_{22}\mu_{11}-\mu_{21}\mu_{12}\right) \sigma S(\ell)\big)\\
  r_2 &= \mu_{12}^{2}\sigma-\mu_{11}^{2}    \\
  \epsilon&=-\mu_{22}\mu_{12}^{2}(\sigma-\bar \sigma),\quad \bar\sigma=\frac{2\mu_{21}\mu_{11}\mu_{12}-\mu_{22}\mu_{11}^{2}}{\mu_{22}\mu_{12}^{2}}.
\end{align*}
Following Section \ref{sec:bezout} it remains to modify $r_1$,
$r_2$ in such a way that $\epsilon$ is replaced by a constant. This
may be done by applying the induction step of
Lemma~\ref{lemma:coprime} once.  To this end, let $\bar r_1,\bar r_2,
\bar p_1,\bar p_2$ be the complex numbers obtained by setting
$\sigma=\bar\sigma$ in the Laplace transforms of $r_1, r_2, p_1,
p_2$. Assume that neither $\bar p_1$ nor $\bar p_2$ are zero. Then the
variables
 \begin{align*}
   q_1&=\frac{\bar p_2 r_1 -\bar r_1 p_2}{\bar p_2\epsilon},\quad
  q_2=\frac{L_{\bar\sigma}(p_1) r_2 -L_{\bar\sigma}(r_2)p_1}{\bar p_1\epsilon}
 \end{align*}
belong to $\Ring_\Rationals$ and, therefore, to $\Ring$. Thus, we have  the \Bezout\ equation $p_1q_1+p_2q_2=1$.

From the above results, one easily verifies that with
$$
 y=q_2\omega_1+q_1\omega_2
$$
one has $\omega_1=p_2y$ and $\omega_2=p_1y$. Hence, $y$ is a basis 
of the system under consideration.

\section{Conclusion}
For a class of convolutional systems associated with boundary coupled second
order partial differential equations we have derived algebraic controllability
results which translate directly into trajectory related controllability
conditions. 
These results rely on a division algorithm for a particular ring of
\Mikusinski\ operators with compact support that is obtained from the operator
solution of the Cauchy problem associated with the given system of partial
differential equations. However, this means that our algebraic setting does not
apply directly to the given boundary value problem but rather to a convolutional
system arising from these solutions in connection with the boundary conditions.
A promising approach allowing an algebraic treatment from the very beginning is
currently under investigation.

The current work was motivated by previous contributions
\cite{BreLoi96msca,Gluesing97siam} in which similar results where
presented for differential delay systems.  Those approaches have been
shown to be useful not only for controllability analysis but also for
the design of closed loop control schemes using the factorization
approach or the method of finite spectrum assignment
\cite{BreLoi97ecc,BreLoi1998mcs,Gluesing2000mcss}. This suggests the
investigation of similar methods for the class of systems considered
within this contribution.

\section*{Acknowledgements}
This paper was written while the first author worked at the
\emph{Laboratoire d'Infor\-ma\-tique} at the \emph{\'Ecole
  Polytechnique} with the financially support of DGA,
\emph{Minist\`ere de la defense francais}. The authors would like to
thank Michel Fliess and Joachim Rudolph for helpful discussions.
\appendix
\renewcommand{\thesection}{\Alph{section}}
\renewcommand{\thesubsection}{\Alph{section}.\arabic{subsection}}
\section{Representation of the operators $S(x)$ and $C(x)$}
In this section we give interpretations of the operators $S(x)$ and
$C(x)$ introduced in Section \ref{bvp_to_sys}. Actually, we restrict
to $S(x)$ from which $C(x)$ can be easily deduced by differentiation
w.r.t.\ $x$. 

If $a> 0$ in equation (\ref{eq:chareq}) we may rewrite $\sigma$ as 
$$
\sigma=\tau^2\left(\left(s+\alpha\right)^2-\beta^2\right),\quad \tau=\sqrt{a},\quad \alpha=\frac{b}{2a},\quad 
\beta=\sqrt{\frac{b^2}{4a^2}-\frac{c}{a}}.
$$
According to \cite{Mik83} the operator $S(x)$ corresponds to the compactly supported function  
\begin{align*}
S(x,t)&=
\left\{\big(h(t+x\tau)-h(t-x\tau)\big)\frac{\mathrm{e}^{-\alpha t}}{2\tau}J_0(\beta\sqrt{\tau^2 x^2-t^2})\right\},
\end{align*}
where $J_0$ denotes the Bessel function of order zero. 
Thus, for any approximate identity $(\varphi_n)$, the operator $S(x)$ possesses a regular representation
given by (cf.\ also App.~\ref{sec:miku_divis})
\begin{align}
 S(x)=f_n(x)/{\varphi_n},\quad f_n(x,t)={\int_{-x\tau}^{x\tau}\frac{\mathrm{e}^{-\alpha \xi}}{2\tau}J_0\big(\beta\sqrt{\tau^2 x^2-\xi^2}\big)\varphi_n(t-\xi)d\xi}.
\end{align}
Contrary, if  $a=0$ in (\ref{eq:chareq}), $S(x)$ can be written as power 
series in the differentiation operator:
$$
 S(x)=\sum_{k=0}^\infty \frac{(s-a)^kx^{2k+1}}{(2k+1)!}.
$$
The convergence of this series is verified directly from the regular representation 
$$
S(x)=\sum_{k=0}^\infty \frac{x^{2k+1} \psi_{k,n}}{(2k+1)!\varphi_{n}},\quad \psi_{k+1,n}=
 \dot\psi_{k,n}-a\psi_{k,n},\quad \psi_{0,n}=\varphi_n
$$
where $(\varphi_{n})$ is any approximate identity of Gevrey order less than $2$. 

\label{sec:app:representation}
\section{Mikusi\'nski Operators and Boehmians}\label{sec:miku}
\subsection{Generalized quotients}\label{sec:quotients}
The set of locally integrable functions with left-bounded support
forms a commutative ring $\mathcal{L}_+$ with respect to the pointwise
addition and the convolution product. A celebrated theorem of
Titchmarsh (\cite[Theorem VII]{Titchmarsh1926}, \cite[Theorem
151]{Titchmarsh1937})\footnote{For several alternative proofs see also
  \cite{Mikusinski1953sm,Mik83,MikBoe87,Yos84}.} states the
following:
\begin{theorem}
  Assume that the convolution product of two locally integrable
  functions $f$ and $g$ the support of which is contained in
  $\Reals^+$ vanishes on $[0,T]$. Then there exist nonnegative real
  numbers $T_g,T_f$ with $T_g+T_f\geqslant T$ such that both, $f$ and $g$,
  vanish identically on $[0,T_f]$ and $[0,T_g]$ respectively.
\end{theorem}
\begin{corollary}\label{cor:titchmarsh}
  The ring $\mathcal{L}_+$ is free of
  divisors of zero.
\end{corollary}
\begin{definition}
  The field $\Mik$ of Mikusi\'nski operators is the quotient field of
  complex-valued locally integrable functions on $\Reals$ with left bounded
  support\footnote{Contrary to this definition $\Mik$ is sometimes defined as the quotient  field of the convolution ring of continuous functions with support in $\Reals^+$.
However, in both cases the obtained fields of fractions are isomorphic.}.
\end{definition}
\begin{definition}
  Consider a commutative ring $R$ together with an $R$-module $M$.
  Let $\Delta$ be a family of sequences of $R$ such that: 
  \begin{enumerate}
    \item  If $(\varphi_n),(\psi_n)\in\Delta$ then $(\varphi_n\psi_n)\in\Delta$.
    \item  For $f,g\in M$ and $(\varphi_n)\in\Delta$  the equality of 
    sequences $(f \varphi_n)=(g \varphi_n)$ implies
     $f=g$.
  \end{enumerate}
  Then the elements of $\Delta$ will be called $\Delta$-sequences in $M$ \cite{Mikusinski1988amh,Mik01}.
\end{definition}
\begin{definition}\label{def:boehmians:general}
Consider an $R$-module $M$ with
  $R$ a commutative ring and $\Delta$ a family of $\Delta$-sequences in $M$. Let
  $\mathscr{A}(M,\Delta)$ the set of all pairs of sequences $(f_n)\in
  M^{\Naturals}$ and $(\varphi_n)\in \Delta$ satisfying $\varphi_i f_j=\varphi_j
  f_i$ for all $i,j\in\Naturals$.  With the (equivalence)-relation $\sim$
  defined by
   $$
    \Big(\big((f_n),(\varphi_n)\big)\sim \big((g_n),(\psi_n)\big)\Big) \Leftrightarrow  
    \Big(\varphi_i g_j=\psi_j f_i \quad\text{for all}\quad i,j\in\Naturals\Big)
   $$
  Then the space $\Boem(M,\Delta)$  of Boehmians on $M$ is defined as 
  $\mathscr{A}(M,\Delta)/\sim$.  For notational simplicity a Boehmian  is simply denoted as $f_n/\varphi_n$.
  The addition on $\Boem(M,\Delta)$ may be defined according to
  $f_n/\varphi_n+g_n/\psi_n=(\psi_n f_n+\varphi_n g_n)/(\varphi_n\psi_n)$ \cite{Mik01,Mikusinski1988amh}.
\end{definition}
In the following $\mathcal{L}(\Reals)$ denotes the space of locally integrable
functions on $\Reals$ and $\mathcal{L}_0(\Reals)$ (resp.
$\mathcal{L}_+(\Reals)$) the subset containing the functions with bounded (resp.\ left-bounded)
support. With the pointwise addition and the convolution product $\mathcal{L}_0$ (resp. $\mathcal{L}_+$)
form a commutative ring and $\mathcal{L}$ (resp.\  $\mathcal{L}_+$) is a $\mathcal{L}_0$-module (resp.\  $\mathcal{L}_+$-module). 

\begin{theoremanddefinition}\cite{Mikusinski1988amh,Mikusinski1990amh}\label{def:boehmians}
Let $\Delta$ be the family of all sequences $(\varphi_n)$ in $\mathcal{L}_0$ satisfying
\begin{enumerate}
  \item  $\exists C\in\Reals: \forall n\in \Naturals:\int_{\Reals}|\varphi_n(t)|dt\leqslant C$  
  \item  $\forall n\in \Naturals:\int_{\Reals}\varphi_n(t) dt=1$ 
  \item  $\forall \varepsilon\in\Reals\,\exists\, n_0\in\Naturals:\forall n\geqslant n_0: \supp {\varphi_n}\subseteq [-\varepsilon,\varepsilon]$.
\end{enumerate}
Then we may define the following spaces of Boehmians:
 \begin{itemize}
 \item The ring $\MikComp$ of \Mikusinski\ operators with compact
   support is defined as 
   $\Boem([\mathcal{L}_0]_{\mathcal{L}_0},\Delta)$
   \cite{Boehme1973ams}.
 \item The ring $\MikReg\supset\MikComp$ of  regular \Mikusinski\ operators  is defined
  to be $\Boem([\mathcal{L}_+]_{\mathcal{L}_+},\Delta)$.
\item The elements of the space
  $\Boem([\mathcal{L}]_{\mathcal{L}_0},\Delta)$, which for short is denoted
  by $\Boem$,  are simply called Boehmians. 
 \end{itemize}
\end{theoremanddefinition}
Obviously, the spaces  $\MikReg$ and $\Boem$ possess  
  the structure of $\MikComp$ modules.
\begin{remark}
  The sets $\MikComp$ and $\MikReg$ regarded as commutative rings are clearly
  isomorphic to subrings of  of $\Mik$. Since the rings $\mathcal{L}_+$ and $\mathcal{L}_0$ are free
  of divisors of zero, for these spaces the equivalence relation in 
  \ref{def:boehmians:general} could be replaced by:    
   $$
    \Big(\big((f_n),(\varphi_n)\big)\sim \big((g_n),(\psi_n)\big)\Big) \Leftrightarrow  
    \Big(\varphi_i g_j=\psi_j f_i \quad\text{for some}\quad i,j\in\Naturals\Big).
   $$
\end{remark}

\subsection{Divisibility  in the ring of \Mikusinski\ operators with compact support}\label{sec:miku_divis}
In this section we shall state some divisibility properties of the ring
$\MikComp$.  For proofs we refer to the cited literature.
\begin{proposition}\label{def:mik:daempfungssatz} 
  The mapping $T^\alpha$ defined by the pointwise multiplication of a function
  $f\in \mathcal{L}_0$ with an exponential function $t\mapsto\mathrm{e}^{\alpha
    t}$ defines an isomorphism on $\mathcal{L}_0$ which extends to an
  isomorphism on $\MikComp$ \cite{Mik83,Woi07diss}.
\end{proposition}
\begin{proposition}\label{def:mik:endwertsatz}
  The mapping $L:\mathcal{L}_{0}\rightarrow\Complexes$ assigning to every
  element of $\mathcal{L}_0$ the value of its integral can be shown to be a
  homomorphism.  Its unique extension to $\MikComp$ is denoted by the same
  symbol.
\end{proposition}
\begin{proposition}\label{thrm:mik:teilbarkeit}
  An operator $a\in\MikComp$ divisible by $(s+\alpha),\,\alpha\in\Complexes$ if
  and only if $L\circ T^{\alpha}(a)=0$.
\end{proposition}
\begin{remark}
  Note that the function $\hat a:\Complexes\rightarrow \Complexes$ given by
  $\hat a(\sigma)=L\circ T^{-\sigma}(a)$ is the Laplace transform $\LaplaceTrans(a)$ of $a$. It is an
  entire function which satisfies a growth condition on the imaginary axis
  derived in \cite{Burzyk1989sm}. Within this context, Proposition
  \ref{thrm:mik:teilbarkeit} means that $s-\alpha$ divides $a$ in $\MikComp$ if
  and only if $\sigma-\alpha$ divides the Laplace transform of $a$ within the
  space of entire functions.
\end{remark}


\bibliographystyle{amsplain}
\bibliography{all}

\end{document}